%
%


\input eplain

\def\parsetimestamp #1: <#2 #3 #4>{{\tt version #2 #3 #4}}

\newif\iftitle
\newtoks\runninghead \runninghead={\hfil}
\newtoks\authors \authors={\hfil}
\newtoks\titlehead
\newtoks\rightheader \newtoks\leftheader \newtoks\pagefoot

\font\hdr=cmssdc10 
\rightheader={\hdr\hfil\the\authors\hfil}
\leftheader={\hdr\hfil\the\runninghead\hfil}
\titlehead={\hfil}
\pagefoot={\hss\tenrm\folio\hss}

\headline={\ifnum\the\pageno<1 \hfil\else
	\iftitle\the\titlehead{\global\titlefalse}\else
		\ifodd\pageno\the\rightheader\else\the\leftheader\fi\fi\fi}
\footline={\ifnum\the\pageno<1 \hfil\else\the\pagefoot\fi}

\fontdimen16\tensy=1\fontdimen17\tensy

\catcode`@ = 11
\newdimen\@StrutBaseDimension
\newdimen\@StrutSkipTemp
\def\ComputeStrut{%
  \@StrutBaseDimension = \baselineskip
  \ifdim\baselineskip < 0pt
    \errhelp = {You probably called \string\offinterlineskip
        before \string\ComputeStrut}
    \errmessage{\string\ComputeStrut: negative
        \string\baselineskip (\the\baselineskip)}%
  \fi}
\def\MyStrut{%
  \vrule height 0.7\@StrutBaseDimension
         depth 0.3\@StrutBaseDimension
         width 0pt}
\ComputeStrut
\catcode`@ = 12

\catcode`\@=11
\def\displaylinesno#1{\displ@y\halign{
	\hbox to \displaywidth{$\@lign\hfil\displaystyle##\hfil$}&
		\llap{$##$}\crcr
	#1\crcr}}
\catcode`\@=12

\newcount\sectno
	\def\newsectno{\global\advance\sectno by1 }

%
%
\outer\def\section#1#2\par{%
	\dimen0=\pagegoal \advance\dimen0 by-120pt
	\removelastskip
	\vskip 20pt plus 40pt \penalty -400 \vskip 0pt plus -32pt
	\newsectno\resno=0\eqnumber=0
	\message{#2}
	\definexref{#1}{\the\sectno}{SECT}	
	\leftline{\bf\the\sectno.\enspace#2}
	\nobreak\bigskip\noindent}	

\outer\def\nonumsection#1\par{
	\resno=0
	\removelastskip
	\vskip 20pt plus 40pt \penalty -400 \vskip 0pt plus -32pt
	\message{#1}
	\leftline{\bf#1}
	\nobreak\bigskip\noindent}

\newskip\resskipamount
	\resskipamount=15pt plus5pt minus3pt
\def\resskip{\vskip\resskipamount}
\def\resbreak{\par\ifdim\lastskip<\resskipamount
	\removelastskip \penalty-300 \resskip\fi}

\def\XXmagic{^^C}
\def\math#1{$#1$}
\def\XX#1 $#2$#3\endXX{\def\XXtest{#2}\ifx\XXtest\XXmagic
	\def\XXnext{#1}\else\def\XXnext{\XX#1\/ \math{#2}#3\endXX}\fi
	\XXnext}

\newcount\resno
	\def\newresno{\global\advance\resno by1 }

\def\result#1#2#3{
	\newresno
	\definexref{#2}{\the\sectno.\the\resno}{#1}
	\resbreak\noindent
	{\bf\the\sectno.\the\resno\enspace{\csname #1word\endcsname.}\enspace}
	{\sl \XX#3 $^^C$\endXX}\par
	\ifdim\lastskip<\bigskipamount%
	\removelastskip\penalty55\bigskip\fi}

\outer\def\theorem#1#2\par{\result{THM}{#1}{#2}}
\outer\def\lemma#1#2\par{\result{LEM}{#1}{#2}}
\outer\def\corollary#1#2\par{\result{COR}{#1}{#2}}

\newtoks\resulttag \resulttag={}

\catcode`@ = 11
\def\Result#1 #2\par{
	\resbreak\noindent
	{\bf\the\sectno.\the\resno\enspace{\csname #1word\endcsname.}\enspace}
	{\sl \XX#2 $^^C$\endXX}\par
	\ifdim\lastskip<\bigskipamount%
	\removelastskip\penalty55\bigskip\fi}

\outer\def\Theorem{\resulttag={THM}\@getoptionalarg\doresult}
\outer\def\Lemma{\resulttag={LEM}\@getoptionalarg\doresult}
\outer\def\Corollary{\resulttag={COR}\@getoptionalarg\doresult}
\def\doresult{%
	\newresno
	\ifx\@optionalarg\empty
	\else
		\definexref{\@optionalarg}{\the\sectno.\the\resno}{\the\resulttag}
	\fi
	\Result{\the\resulttag}
}
\catcode`@ = 12

\def\proof{\noindent{{\sl Proof. }}}

\def\sqr#1#2{{\vbox{\hrule height.#2pt
    \hbox{\vrule width.#2pt height#1pt \kern#1pt
        \vrule width.#2pt}\hrule height.#2pt}}}
\def\eqed{\sqr53}
\def\qed{%
	\ifmmode\eqno\eqed
	\else\nobreak\ \hfill\eqed\medbreak\fi}

%
%


\input amssym.def
\input amssym

\font\ninerm=cmr9 \font\ninesy=cmsy9
\font\ninei=cmmi9
\font\nineit=cmti9 \font\ninebf=cmbx9
\font\ninesl=cmsl9
\def\ninepoint{\textfont0=\ninerm\textfont1=\ninei
  \textfont2=\ninesy\scriptfont0=\sevenrm\scriptfont1=\seveni
  \scriptfont2=\sevensy\scriptscriptfont0=\fiverm\scriptscriptfont1=\fivei
  \scriptscriptfont2=\fivesy\def\rm{\fam0\ninerm}\let\it=\nineit\let\bf=\ninebf
  \let\sl=\ninesl\baselineskip=13pt\rm}

\def\al{\alpha}
\def\be{\beta}
\def\de{\delta}
\def\De{\Delta}

\def\ga{\gamma}

\def\la{\lambda}
\def\La{\Lambda}

\def\sg{\sigma}

\def\Th{\Theta}

%
%

\def\cN{{\cal N}}

\def\cV{{\cal V}}

%
%

\def\cx{{\Bbb C}}
\def\fld{{\Bbb F}}


%
%

\def\inv{^{-1}}
\def\sbs{\subseteq}

\def\seq#1#2#3{#1_{#2},\ldots,#1_{#3}}

%
%
\def\ip#1,#2{\langle#1,#2\rangle}
\def\one{{\bf1}}

\def\tr{\mathop{\hbox{\rm tr}}\nolimits}

%
%


\def\sym#1{{{\rm Sym}(#1)}}

\magnification=\magstep1
\hoffset=1.2cm
\hsize=27pc
\vsize=42pc
\baselineskip=12pt
\normalbaselineskip=12pt

\authors={Chan, Godsil and Munemasa}
\runninghead={Jones Pairs}

\def\mnf{M_k({\fld})}
\def\del#1{\De_{#1}}
\def\deb{\del B}
\def\dxd#1#2#3{\del{#1}X_{#2}\del{#3}}
\def\xdx#1#2#3{X_{#1}\del{#2}X_{#3}}

\def\eij{E_{i,j}}
\def\nom#1{\cN_{#1}}
\def\nab{\nom{A,B}}
\def\nabp{\nab'}
\def\na{\nom A}
\def\sgip{\sg_{i+1}}
\def\snv{^{(-)}}
\def\snt{^{(-)T}}

\def\thab{\Th_{A,B}}
\def\vom{V^{\otimes m}}

\def\xa{X_A}
\def\yij{Ae_i\circ Be_j}

\def\Mnf{M_k(\fld)}

\def\eqconstruct#1{\the\sectno.#1}

\titletrue
\null\vskip1.0cm
\centerline{\bf Four-Weight Spin Models and Jones Pairs}
\vskip1.0cm
\centerline{ 
Ada Chan,
Chris Godsil
\footnote{$^1$}{Support from a National Sciences and Engineering 
Council of Canada operating grant is gratefully acknowledged.}
and Akihiro Munemasa}

\footnote{}{2000 Mathematics Subject Classification 05E30}

\vskip0.5cm
\centerline{Combinatorics and Optimization} 
\centerline{University of Waterloo}
\centerline{Waterloo, Ontario}
\centerline{Canada N2L 3G1}

\vskip1.0cm
\centerline{ABSTRACT}
\bigskip
{\narrower\noindent 
We introduce and discuss Jones pairs.  These provide a generalization
and a new approach to the four-weight spin models of Bannai and
Bannai.  We show that each four-weight spin model determines a ``dual''
pair of association schemes.
\par}

\section{pairs}
Jones Pairs

The space of $k\times k$ matrices acts on itself in three distinct
ways: if $C\in\mnf$, we can define endomorphisms $X_C$, $\De_C$ and
$Y_C$ by
$$
X_C(M) =CM,\qquad \De_C(M)=C\circ M, \qquad Y_C(M)=MC^T.
$$
If $A$ and $B$ are $k\times k$ matrices, we say $(A,B)$ is a {\sl
one-sided Jones pair}\/ if $\xa$ and $\deb$ are invertible and
$$
\xdx ABA = \dxd BAB.\eqdef{braid}
$$
We call this the {\sl braid relation}, for reasons that will become
clear as we proceed.  We note that $\xa$ is invertible if and only if
$A$ is, and $\deb$ is invertible if and only if the Schur inverse
$B\snv$ is defined.  (Recall that if $B$ and $C$ are matrices of the
same order, then their Schur product $B\circ C$ is defined by the
condition
$$
(B\circ C)_{i,j}=B_{i,j}C_{i,j}
$$
and $B\circ B\snv=J$.)  We will see that each one-sided Jones pair
determines a representation of the braid group $B_3$, and this is one
of the reasons we are interested in Jones pairs.  The pair $(I,J)$
forms a trivial but useful example.  A one-sided Jones pair $(A,B)$ is
{\sl invertible} if $A\snv$ and $B\inv$ both exist.

We observe that $X_A$ and $Y_A$ commute and that
$$
Y_A\deb Y_A =\deb Y_A\deb\eqdef{braid2}
$$
if and only if $(A,B^T)$ is a Jones pair.  A pair of matrices $(A,B)$
is a {\sl Jones pair} if both $(A,B)$ and $(A,B^T)$ are one-sided
Jones pairs.  From a Jones pair we obtain a family of representations
of the braid groups $B_r$, for all $r$.

In this paper we describe the basic theory of Jones pairs.  We find
that if $A\snv$ exists, then $(A,A\snv)$ is a Jones pair if and only
if $A$ is a spin model in the sense of Jones \cite{MR89m:57005}, and
that invertible Jones pairs correspond to the 4-weight spin models due
to Bannai and Bannai \cite{MR97f:57005}.  The theory we develop
includes the basic theory of these spin models.  Finally we prove that
each invertible Jones pair $(A,B)$ of $n\times n$ matrices determines
a dual pair of association schemes, one built from $A$ and the other
from $B$.

\section{brdreps}
Representations of the Braid Group

The braid group $B_n$ on $n$ strands has $n-1$ generators
$\seq{\sg}1{n-1}$ which satisfy the relations
$$
\sg_i\sgip\sg_i =\sgip\sg_i\sgip
$$
and, if $|i-j|>1$,
$$
\sg_i\sg_j =\sg_j\sg_i.
$$
We will show how we can use Jones pairs to construct a class of
representations of the braid group and, in certain cases, obtain
invariants of oriented links.  These constructions are due to Jones
\cite{MR89m:57005}, who did not feel the need to write out proofs.
We did, and so we record them here.  Nonetheless, these proofs play no
role in the sections that follow (and there will be no exam).

Let $\vom$ denote the tensor product of $m$ copies of the
$n$-dimensional vector space $V$.  Let $(A,B)$ be a pair of $n\times
n$ matrices.  Let $e_i$ denote the $i$-th standard basis vector of
$\fld^n$ and let $\eij$ be the matrix $e_ie_j^T$.  Let $g_{2i-1}$ be
the element of ${\rm End}(\vom)$ we get by applying $A$ to the $i$-th
tensor factor of $\vom$.  We let $g_{2i}$ denote the element of ${\rm
End}(\vom)$ such that
$$
g_{2i}(e_{r_1}\otimes\cdots\otimes e_{r_m})
	=B_{r_i,r_{i+1}}(e_{r_1}\otimes\cdots\otimes e_{r_m}).
$$

We remark that $(A,B)$ is a one-sided Jones pair if and only if the
map sending $\sg_1$ to $g_1$ and $\sg_2$ to $g_2$ determines a
representation of $B_3$.  Our next result is from Jones
\cite[Section~3.3]{MR89m:57005}.

\lemma{brd-rep}
Let $(A,B)$ be a pair of $n\times n$ matrices and let the
endomorphisms $g_i$ be defined as above.  If $r\ge4$ and
$m\ge\lfloor(r+1)/2\rfloor$, the following are equivalent.
\smallskip
\item{(i)} $(A,B)$ is a Jones pair.
\item{(ii)} The mapping that assigns $\sg_i$ in $B_r$ to $g_i$
defines a representation of $B_r$ in ${\rm End}(\vom)$.\qed

\proof
It suffices to show that, for each positive integer $r\geq4$, 
we have
$$
g_{i}g_{i+1}g_{i}
	=g_{i+1}g_{i}g_{i+1}\quad(i=1,\ldots,r-2).\eqdef{eq:braidg}
$$
if and only if \eqref{braid} and \eqref{braid2} hold.
By the
definition of the action of $g_i$, we see that the equations
\eqref{eq:braidg} are equivalent to
$$
g_{1}g_{2}g_{1}=g_{2}g_{1}g_{1},\quad
g_{2}g_{3}g_{2}=g_{3}g_{2}g_{3}\eqdef{eq:brdg4}
$$
on $V\otimes V$. If we identify $V\otimes V$ with $\Mnf$ via the
map $\phi:e_{i}\otimes e_{j}\mapsto E_{ij}$, then
$$\eqalign{
\phi g_{1}(e_{i}\otimes e_{j})&=\phi(Ae_{i}\otimes e_{j})\cr
	&=\phi(\sum_{h=1}^{n}A_{hi}e_{h}\otimes e_{j})\cr
	&=\sum_{h=1}^{n}A_{hi}E_{hj}\cr
	&=AE_{ij}\cr
	&=\xa(E_{ij}),\cr}
$$

$$\eqalign{
\phi g_{2}(e_{i}\otimes e_{j})&=\phi(B_{ij}e_{i}\otimes e_{j})\cr
	&=B_{ij}E_{ij}\cr
	&=\deb(E_{ij}),\cr}
$$

$$\eqalign{
\phi g_{3}(e_{i}\otimes e_{j}) &=\phi(e_{i}\otimes Ae_{j})\cr
	&=\phi(\sum_{h=1}^{n}A_{hj}e_{i}\otimes e_{h})\cr
	&=\sum_{h=1}^{n}A_{hj}E_{ih}\cr
	&=\sum_{h=1}^{n}A^{T}_{jh}E_{ih}\cr
	&=E_{ij}A^{T}\cr
	&=Y_A(E_{ij}).\cr}
$$
Therefore, we have shown
$$
g_{1}=\phi^{-1}\xa\phi,\quad
g_{2}=\phi^{-1}\deb\phi,\quad
g_{3}=\phi^{-1}Y_A\phi,
$$
and hence the equations \eqref{eq:brdg4} are equivalent to
\eqref{braid} and \eqref{braid2}.\qed

Although representations of the braid group are interesting in their
own right, the work we have discussed is motivated by the relation
with link invariants.  The following lemma gives a sufficient
condition for a Jones pair to yield a link invariant.  (We comment on
the constraints on $A$ and $B$ following the proof of \ref{binvasnv}.)

\lemma{markov}
Let $(A,B)$ be a Jones pair of $n\times n$ matrices such that
for $i=1,\ldots,n$ we have
$$
A_{i,i}=(A\inv)_{i,i}=1/\sqrt{n}
$$
and
$$
\sum_{j=1}^n B_{ij} =\sum_{j=1}^n B_{ij}\inv=\sqrt{n}.
$$
If $\seq g1{r-1}$ are defined as above and 
$h\in\langle\seq g1{r-2}\rangle$, then
$$
\tr(hg_{r-1}^{\pm1})={1\over n}\tr(h)\tr(A)^{\pm1}.
$$

\proof
Let $r'=\lfloor(r+1)/2\rfloor$ and let $\cV$ denote the tensor product
of $r'$ copies of $V$.  Suppose first that $r$ is even.  If
$\alpha\in\{1,\ldots,n\}^{r'-1}$, denote by $e_{\alpha}$ the vector
$e_{\alpha_{1}}\otimes\cdots\otimes e_{\alpha_{r'-1}}$. Then
$$
\{e_{\alpha}\otimes e_{i} \mid \alpha\in\{1,\ldots,n\}^{r'-1},\;
1\leq i \leq n \}
$$
forms a basis of $\cV$. Since $h$ leaves the subspace
$$
\langle
e_{\alpha}\otimes e_{i} \mid \alpha\in\{1,\ldots,n\}^{r'-1}
\rangle
$$
invariant, we can write
$$
h(e_{\alpha}\otimes e_{i})=\sum_{\beta} h_{\beta,\alpha}^{i}
e_{\beta}\otimes e_{i},
$$
where $h_{\beta,\alpha}^{i}\in\fld$. Now we have
$$
\tr(h)=\sum_{i=1}^{n}\sum_{\alpha} h_{\alpha,\alpha}^{i}.
$$
Since
$$\eqalign{
hg_{r-1}^{\pm1}(e_{\alpha}\otimes e_{i})
	&=h(e_{\alpha}\otimes \sum_{j=1}^{n} (A^{\pm1})_{ji}e_{j})\cr
	&=\sum_{j=1}^{n} (A^{\pm1})_{ji} \sum_{\beta} h_{\beta,\alpha}^{j}
		e_{\beta}\otimes e_{j},\cr}
$$
we have
$$\eqalign{
\tr(hg_{r-1}^{\pm1})&=\sum_{i=1}^{n} \sum_{\alpha} (A^{\pm1})_{ii}
	h_{\alpha,\alpha}^{i}\cr 
	&={1\over\sqrt{n}}
		\sum_{i=1}^{n} \sum_{\alpha} h_{\alpha,\alpha}^{i}\cr
	&={1\over n} \tr(h)\tr(A).\cr}
$$
This proves the lemma when $r$ is even.

Next suppose $r$ is odd.  For $\alpha\in\{1,\ldots,n\}^{r'-2}$, denote
by $e_{\alpha}$ the vector $e_{\alpha_{1}}\otimes\cdots\otimes
e_{\alpha_{r'-2}}$.  Then
$$
\{e_{\alpha}\otimes e_{i}\otimes e_{j} \mid
	\alpha\in\{1,\ldots,n\}^{r'-2},\; 1\leq i,j \leq n \}
$$
forms a basis of $\cV$.  Since $h$ acts as the identity on the last
tensor factor of $\cV$, we can write
$$
h(e_{\alpha}\otimes e_{i}\otimes e_{j})=
\sum_{k=1}^{n}\sum_{\beta} h_{(\beta,k),(\alpha,i)}
e_{\beta}\otimes e_{k}\otimes e_{j},
$$
where $h_{(\beta,k),(\alpha,i)}\in\fld$ are independent of $j$.
Now we have
$$
\tr(h)=n\sum_{i=1}^{n}\sum_{\alpha} h_{(\alpha,i),(\alpha,i)}.
$$
Since
$$\eqalign{
hg_{r-1}^{\pm1}(e_{\alpha}\otimes e_{i}\otimes e_{j})
	&=B_{ij}^{\pm1} h(e_{\alpha}\otimes e_{i}\otimes e_{j})\cr
	&=B_{ij}^{\pm1} \sum_{k=1}^{n} \sum_{\beta} h_{(\beta,k),(\alpha,i)}
		e_{\beta}\otimes e_{k}\otimes e_{j},\cr}
$$
we have
$$\eqalign{
\tr(hg_{r-1}^{\pm1}) &=\sum_{i,j=1}^{n} \sum_{\alpha} B_{ij}^{\pm1} 
		h_{(\alpha,i),(\alpha,i)}\cr
	&=\sum_{i=1}^{n} \sum_{\alpha} h_{(\alpha,i),(\alpha,i)}
		\sum_{j=1}^{n} B_{ij}^{\pm1}.\cr}
$$
Given the constraints on $B$ and $B\snv$, we see that the lemma holds
when $r$ is odd.\qed

\section{more}
Further Properties

We develop some basic properties of one-sided Jones pairs.  We begin
with an alternative form of the definition.

Equation \eqref{braid} is equivalent to the condition that, for all
matrices $M$ in $\mnf$,
$$
A(B\circ (AM)) =B\circ(A(B\circ M)).\eqdef{brd2}
$$

If $M$ and $N$ are elements of $\mnf$, then $\tr(M^TN)$ is a
non-degenerate bilinear form on $\mnf$.  If $Y\in{\rm End}(\mnf)$,
denote its adjoint relative to this form by $Y^T$, and call it the
transpose of $Y$.  It is easy to
verify that 
$$
\xa^T =X_{A^T},\qquad \deb^T =\deb.
$$

\lemma{jp-equivs}
If $(A,B)$ is a one-sided Jones pair, so is:
\item{(a)} $(A^T,B)$.
\item{(b)} $(A\inv,B\snv)$.
\item{(c)} $(D\inv AD,B)$, where $D$ is invertible and diagonal
\item{(d)} $(A,BP)$, where $P$ is a permutation matrix.
\item{(e)} $(P\inv AP,P\inv BP)$, where $P$ is a permutation matrix. 
\item{(f)} $(\la A,\la B)$, for any non-zero $\la$ in $\fld$.

\proof
The first claim follows by taking the transpose of Equation
\eqref{braid}, and the second by taking the inverse and noting that
$$
\xa\inv =X_{A\inv},\qquad \deb\inv =\del{B\snv}.
$$
If $D$ is diagonal and invertible, then
$$
X_{D\inv}X_C X_D =X_{D\inv CD},\qquad X_{D\inv}\De_C X_D=\De_C.
$$
Then (c) follows by conjugating \eqref{braid} above by $X_D$.

Next, if $P$ is a permutation matrix, then
$$
(BP)\circ M =(B\circ MP\inv)P
$$
and 
$$
(B\circ C)P= BP \circ CP,
$$
hence
$$\eqalignno{
A(BP\circ(AM)) &=A(B\circ(AMP\inv))P\cr
	&=(B\circ A(B\circ(MP\inv)))P\cr
	&=(BP\circ A(B\circ(MP\inv)P))\cr
	&=(BP\circ A(BP\circ(M))\cr}
$$
We leave (e) as an exercise, while (f) is truly trivial.\qed

\section{evecs}
Eigenvectors

Note that we can rewrite the braid relation in the equivalent form
$$
\deb\inv\xa\deb =\xa\deb\xa\inv,
$$
from which we see that the endomorphisms $\xa$ and $\deb$ are similar.
Since
$$
\deb\eij =B_{i,j}\eij,
$$
$\eij$ is an eigenvector for the operator $\deb$, with eigenvalue $B_{i,j}$.
Since $\xa$ and $\deb$ are similar, we conclude that $\xa$ is
diagonalizable, and its eigenvalues are the entries of $B$.

\lemma{evecs}
The matrices $A$ and $B$ form a one-sided Jones pair if and only if,
for all $i$ and $j$ we have 
$$
A(Ae_i\circ Be_j) =B_{i,j}(Ae_i\circ Be_j).
$$

\proof
Put $M=\eij$ in the second form, \eqref{brd2}, of the braid
relation.  This yields
$$
A(B\circ (A\eij)) =B\circ(A(B\circ\eij)).\eqdef{abaeij}
$$
Observe that
$$
B\circ(A\eij) =B\circ(Ae_ie_j^T) =(Ae_i\circ Be_j)e_j^T
$$
and $B\circ\eij =B_{i,j}\eij$.  Therefore,
$$
B\circ(A(B\circ\eij)) =B_{i,j} B\circ(A\eij) 
	=B_{i,j}(Ae_i\circ Be_j)e_j^T
$$
and hence the lemma follows from \eqref{abaeij}.\qed

\lemma{bcols}
If $(A,B)$ is a one-sided Jones pair, then $B$ has constant column sum.

\proof
We have
$$
A(\yij) =B_{i,j}(\yij).
$$
Since $B\snv$ exists and $A$ is invertible, the vectors
$$
Ae_r\circ Be_j,\quad r=1,\ldots,n
$$
are linearly independent, and therefore they form a basis of $\fld^n$
consisting of eigenvectors of $A$.  Consequently, the entries $B_{r,j}$
in the $j$-th column of $B$ provide a complete list of the eigenvalues
of $A$.  Thus each column of $B$ sums to $\tr(A)$.\qed

It follows from this lemma that if $(A,B^T)$ is a Jones pair, then the
row sums of $B$ are also constant.  The argument used will also show
that, if $(A,B)$ is a one-sided Jones pair and $A\snv$ and $B\inv$
exist, then the row sums of $B$ are constant.

\medbreak
\lemma{ada}
Let $(A,B)$ be a pair of $n\times n$ matrices and let $D_j$ be the
diagonal matrix whose $r$-th diagonal entry is $B_{r,j}$.  Then
$(A,B)$ is a one-sided Jones pair if and only if
$$
AD_jA =D_jAD_j
$$
for $j=1,\ldots,n$.

\proof
By \ref{evecs} we have that $(A,B)$ is a one-sided Jones pair if and
only if
$$
A\,\yij =B_{i,j}\,\yij
$$
for all $i$ and $j$.
Since $Be_j\circ x= D_jx$, this is equivalent to
$$
AD_j\, Ae_i =B_{i,j} D_j Ae_i,
$$
from which the lemma follows.\qed

A representation of $B_3$ over $\fld$ is given by assigning invertible
$n\times n$ matrices $A_1$ and $A_2$ over $\fld$ to the generators
$\sg_1$ and $\sg_2$ of $B_3$, such that the braid condition holds:
$$
A_1A_2A_1 =A_2A_1A_2.
$$
It follows from \ref{ada} that if $(A,B)$ is a one-sided Jones pair
and $D_j$ is the diagonal matrix whose $r$-th diagonal entry is
$B_{r,j}$, then the pair of matrices $(A,D_j)$ provide a representation
of $B_3$.

Conversely, suppose $A$ and $D$ are invertible matrices that provide a
representation of $B_3$, and that $D$ is diagonal.  Set $B=DJ$.
Then $(DJ)\circ M=DM$ for any matrix $M$, and therefore
$$
A(B\circ(AM)) =A(DJ\circ(AM)) =ADAM =DADM =B\circ(A(B\circ M)).
$$
Thus $(A,D)$ determines a representation of $B_3$ if and only if
$(A,DJ)$ is a one-sided Jones pair.  To get an invertible one-sided
Jones pair containing $A$, we require $n$ linearly independent
diagonal matrices $D_i$ such that $(A,D_i)$ gives a representation of
$B_3$ for $i=1,\ldots,n$.

\section{exch}
The Exchange Lemma
 
It is trivially true that, if $Ae_i\circ Be_j$ is an eigenvector for
$A$, then so is $Be_j\circ Ae_i$.  Somewhat surprisingly, this yields a
very useful ``Exchange Lemma''.

\lemma{xch}
If $A$, $B$, $C$ and $Q$, $R$, $S$ are elements of $\mnf$, then
$$
\xdx ABC =\dxd QRS
$$
if and only if 
$$
\xdx ACB =\dxd RQ{S^T}.
$$

\proof
The first relation holds if and only if 
$$
A(B\circ(C\eij) =Q\circ(R(S\circ\eij))\eqdef{abcqrs}
$$
We have
$$
B\circ(C\eij) =(Ce_i\circ Be_j)e_j^T
$$
and therefore
$$
A(B\circ(C\eij)) =A(Ce_i\circ Be_j)e_j^T
$$
Further,
$$
Q\circ(R(S\circ\eij)) = S_{i,j}(Re_i\circ Qe_j)e_j^T.
$$
Consequently, our first relation is equivalent to the system of
equations
$$
A(Ce_i\circ Be_j) =S_{i,j}(Re_i\circ Qe_j).
$$

To complete the argument, note that
$$
A(Ce_i\circ Be_j)e_i^T =A(C\circ(BE_{j,i}))
$$
and
$$
S_{i,j}(Re_i\circ Qe_j)e_i^T =S_{i,j}(R\circ QE_{j,i})
	=R\circ Q(S^T\circ E_{j,i}).
$$
We conclude that \eqref{abcqrs} holds if and only if 
$$
A(C\circ(BE_{j,i})) =R\circ Q(S^T\circ E_{j,i}),
$$
which is equivalent to the second of the two relations in the
lemma.\qed 

We note one corollary of this lemma.

\corollary{a3b3}
Let $A$ and $B$ be matrices such that $A\inv$ and $B\snv$ exist.  Then
$(A,B)$ is a one-sided Jones pair if and only if
$$
\xdx AAB =\dxd AB{B^T}.
$$

\proof
Apply the exchange lemma to the braid relation 
$\xdx ABA =\dxd BAB$.\qed

This immediately yields the following:

\lemma{3reg}
Let $A$ and $B$ be invertible and Schur invertible matrices.  Then
$(A,B)$ is a one-sided Jones pair if and only if
$$
\dxd{A\snv}AA =\xdx B{B^T}{B\inv}.
$$

\section{dual}
Duality

Let $A$ and $B$ be two $n\times n$ matrices.  We define $\nab$ to be
the space of $n\times n$ matrices such that $\yij$ is an eigenvector
for all $i$ and $j$; this is an algebra under matrix multiplication.
We call it the {\sl Nomura algebra} of the pair $(A,B)$.  If
$R\in\nab$ and $S$ is the $n\times n$ matrix such that
$$
R\,\yij= S_{i,j}\,\yij,
$$
then we denote $S$ by $\thab(R)$.  We denote the image of $\nab$ under
$\thab$ by $\nabp$, and observe that this is a commutative algebra
under Schur multiplication.  In all cases, $I\in\nab$ and
$J=\thab(I)\in\nabp$.  If $(A,B)$ is a one-sided Jones pair, then
$\nab$ contains all polynomials in $A$ and $\nabp$ contains all
``Schur polynomials'' in $B$.

\lemma{}
Let $(A,B)$ be a one-sided Jones pair and let $\La$ denote the
operator $\xdx ABA=\dxd BAB$.  Then $\La^2$ commutes with $\xa$
and $\deb$.  Further $\La\inv\xa\La=\deb$ and $\La\inv\deb\La=\xa$.

\proof
For the first claim, note that
$$
\La^2 =(\xa\deb)^3=(\xdx ABA)^2,
$$
whence $\La^2$ commutes with $\xa\deb$ and $\xa\deb\xa$.  So $\La^2$
commutes with $\xa$ and $\deb$.  Next
$$
\xa\La =\xa\dxd BAB =\La\deb
$$
and
$$
\deb\La=\deb\xdx ABA =\La\xa.\qed
$$

We cannot prove that conjugation by $\La$ swaps $\nab$ and $\nabp$,
but the following is a very useful consolation.

\theorem{nabchar}
Let $(A,B)$ be a pair of $n\times n$ matrices.  Then $R\in\nab$ and
$\thab(R)=S$ if and only if $X_R\deb\xa=\deb\xa\del{S}$.  If $A$ is
invertible and $B$ is Schur invertible, then $\thab$ is an
isomorphism.

\proof
We have
$$
X_R\deb\xa(E_{i,j})= R\,(\yij)\,e_j^T
$$
and
$$
\deb\xa\del{S}(E_{i,j}) =S_{i,j}\,(\yij)\,e_j^T,
$$
from which the result follows the first claim follows.

Since $\deb$ and $\del{S}$ commute we see that
$$
X_R\deb\xa=\deb\xa\del{S}
$$
if and only if
$$
X_R\deb\xa\deb=\deb\xa\deb\del{S},
$$
i.e., if and only if $X_R\La=\La\del{S}$.  If $A\inv$ and $B\snv$
exist, then $\La$ is invertible, whence it follows that $\thab$ is an
isomorphism. \qed

Since $\nabp$ is commutative, we see that $\nab$ is a commutative
algebra.  We list three equivalent forms of the first part of this
theorem, for later use.

\corollary{rbabas}
If $R\in\nab$ and $\thab(R)=S$, then 
\smallskip
\item{(a)} $\xdx {B^T}A{R^T}=\dxd{S^T}{B^T}{A}$. 
\item{(b)} $\dxd{R}{B\snv}{B^T} =\xdx A{A\inv}S$.
\item{(c)} $\dxd{B^T}{B\snt}{R} =\xdx{S^T}{A\inv}{A^T}$.

\proof
We have
$$
\xdx RBA =\dxd BAS,\eqdef{rba1}
$$
and so using the exchange identity,
$$
\xdx RAB =\dxd AB{S^T}.
$$
Taking the transpose of each side, we get
$$
\xdx{B^T}A{R^T} =\dxd{S^T}{B^T}A,
$$
which yields (a).

Next rewrite \eqref{rba1} as
$$
\dxd{B\snv}RB =\xdx A{S}{A\inv}
$$
and apply the exchange lemma to get (b).  Taking the
transpose of each side yields (c).\qed

\section{type2}
Type II Matrices

An $n\times n$ matrix $A$ is a {\sl type II matrix}\/ if $A\snv$
exists and
$$
A A\snt =nI.
$$
Hadamard matrices provide one important class of examples.  When
discussing type II matrices we will assume implicitly that $n$ is
coprime to the characteristic of our underlying field $\fld$.  The
algebra $\nom{A,A\snv}$ is known as the {\sl Nomura algebra} of $A$.
We will usually denote it by $\na$ (rather than $\nom{A,A\snv}$, and
similarly we abbreviate $\Th_{A,A\snv}$ to $\Th_A$.

We recall that a one-sided Jones pair $(A,B)$ is {\sl invertible} if
$A\snv$ and $B\inv$ both exist.  Somewhat surprisingly, invertibility
implies that $A$ and $B$ are type II matrices.

\theorem{jp-type2}
Suppose $(A,B)$ is a one-sided Jones pair.  If $B\inv$ exists, then
the diagonal of $A$ is constant and $A$ and $B$ are type II matrices.

\proof
We use \ref{rbabas} which implies that, since $A\in\nab$ and
$\thab(A)=B$, 
$$
\dxd {B^T}{B\snt}A =\xdx {B^T}{A\inv}{A^T}.
$$
Apply both sides to $I$.  On the left we get
$$\eqalign{
\dxd {B^T}{B\snt}A (I) &=B^T\circ(B\snt(A\circ I))\cr
	&=(B^T\circ B\snt)(A\circ I)\cr
	& =J(A\circ I).\cr}
$$
and, on the right,
$$
\xdx {B^T}{A\inv}{A^T}(I)=B^T(A\inv\circ A^T).
$$
Therefore,
$$
J(A\circ I) =B^T(A\inv\circ A^T).
$$
Since $(A,B)$ is a one-sided Jones pair, $B^T\one=\be\one$, where
$\be=\tr(A)$.  If $B$ is invertible, then $\be\ne0$ and we have
$$
A\inv\circ A^T =(B^T)\inv J(A\circ I) =\be\inv J(A\circ I).
$$
The sum of the entries in the $i$-th column of $A\inv\circ A^T$ is
$$
\sum_r (A\inv)_{r,i}(A^T)_{r,i} =\sum_r(A\inv)_{r,i}A_{i,r}=1,
$$
from which it follows that all columns of $J(A\circ I)$ are equal.
Therefore, the diagonal of $A$ is constant, and so $A\inv\circ A^T$ is
a multiple of $J$.  Since the columns of $A\inv\circ A^T$ sum to 1, we
conclude that $nA\inv\circ A^T=J$, and therefore $A$ is a type II
matrix.

Now, we prove that $B$ is type II.  We know that $A$ and $B$ are both
invertible and Schur invertible so, from \ref{3reg},
$$
\dxd{A\snv}AA =\xdx B{B^T}{B\inv}
$$
Applying each side to $I$, we find that
$$
A\snv\circ(A(A\circ I))= B(B^T\circ B\inv)
$$
As the diagonal of $A$ is constant, this implies that $B(B^T\circ
B\inv)$ is a multiple of $J$, hence
$$
B^T\circ B\inv =cB\inv J,
$$
for some $c$.  Here we can write the right side as $DJ$, where $D$ is
diagonal.  However, arguing as before, all rows of $B^T\circ B\inv$
sum to 1, and it follows that $B$ is type II.\qed

\lemma{binvasnv}
If $(A,B)$ is a Jones pair and $A\snv$ exists, then $B\inv$ exists
(and $A$ and $B$ are type II matrices).

\proof
We use \ref{rbabas}(b) with $A=R$ and $B=S$ to get
$$
\del{A}X_{B\snv}\del{B^T} =\xa\del{A\inv}X_B.
$$
Applying each side of this to $J$ yields
$$
A\circ(B\snv B^T) =A(A\inv\circ(BJ)).
$$
Since $(A,B)$ is a Jones pair, the row sums of $B$ are all equal, so
the right side here is a multiple of $I$.  Since the diagonal of
$B\snv B^T$ is constant and $A\snv$ exists, $B\snv B^T$ must be a scalar
matrix.  Therefore, $B$ is invertible (and type II).\qed

It follows that if $(A,B)$ is an invertible Jones pair, then all
rows and all columns of $B$ have the same sum.  From the proof of
\ref{bcols}, this sum is $\tr(A)$ and, from \ref{jp-type2}, if $B$ is
invertible, then the diagonal of $A$ is constant.  From
\ref{jp-equivs}(b), we know that $(A\inv,B\snv)$ is a one-sided
Jones pair if $(A,B)$ is.  Therefore, by \ref{jp-equivs}(b), we
conclude that if $(A,B)$ is an invertible Jones pair then there is a
non-zero scalar $\la$ such that the pair $(\la A,\la B)$ satisfies the
conditions of \ref{markov}.

The next result is due to Jaeger, Matsumoto and Nomura
\cite{MR99f:05125}.  We include a short new proof of it here, using the
exchange lemma.

\lemma{nancor}
Let $A$ be a type II matrix.  Then $R\in\na$ if and only if
$\Th_A(R)\in\nom{A^T}$.   Further, if $R\in\na$, then
$\Th_{A^T}(\Th_A(R))=nR^T$. 

\proof
Applying \ref{rbabas}(b) (with $B=A\snv$) we have that
$$
\xdx{A}{A\inv}{S}=\dxd{R}{A}{A\snt}
$$
and therefore
$$
\xdx S{A^T}{A\inv} =\dxd {A^{-1(-)}}{A\inv}R.\eqdef{sa4r}
$$
Since $A$ is type II we have $nA\inv=A\snt$, and so \eqref{sa4r}
yields that $S\in\nom{A\inv}$ and $\Th_{A\inv}(S)=nR$.  Finally,
$\nom{A\inv,A\inv{}\snv}=\nom{A\snt,A^T}$ and therefore
$\Th_{A\inv}(S)=\Th_{A^T}(S)^T$.\qed

Since $A^T$ is type II if $A$ is, this lemma implies that $\Th_{A^T}$
maps $\nom{A^T}$ into $\na$.  If we apply the exchange lemma to the
transpose of \eqref{sa4r}, we see that $R^T\in\na$ and
$\Th_A(R^T)=\Th_A(R)^T$. 
Therefore $\na$ and $\nom{A^T}$ are closed under transposes.  Since
$$
\Th_A(Q)\circ\Th_A(R)=\Th_A(QR),
$$
we also find that $\nom{A^T}$ is closed under the Schur product, as
well as under multiplication.  Since $I$ and $J$ both lie in
$\nom{A^T}$, it follows that $\nom{A^T}$ is the Bose-Mesner algebra of
an association scheme.  Since $\na$ is the image of $\nom{A^T}$ under
$\Th_{A^T}$, it too is Schur-closed and is the Bose-Mesner algebra of
a second association scheme.  (These schemes necessarily form a dual
pair, but we do not stop to discuss this.)

\lemma{nab-jpr}
Let $(A,B)$ be a pair of type II matrices of the same order.  Suppose
the diagonal of $A$ is constant, and all row sums of $B$ are equal.
If $A\in\nab$, then there is a scalar $c$ such that $(A,cB)$ is a
one-sided Jones pair.

\proof
If $A\inv$ exists and $A\in\nab$, then $A\inv\in\nab$.  
Hence there is a matrix $S$ such that
$$
\xdx{A\inv}BA =\dxd BAS
$$
and from the exchange lemma it follows that
$$
\xdx{A\inv}AB=\dxd AB{S^T}.
$$
If we apply both sides of this equality to $J$ we find that
$$
A\inv(A\circ(BJ)) =A\circ(BS^T)
$$

Suppose $BJ=\be J$.  Then the left side above is equal to $\be I$, and
therefore $\be A\snv\circ I=BS^T$.  If $A\circ I=\al I$, this shows
that
$$
\thab(A\inv) ={\be\over\al} B^{-T} ={\be\over n\al}B\snv
$$
and consequently
$$
\thab(A) ={n\al\over\be}B.
$$
This implies that
$$
\xdx ABA ={n\al\over\be}\dxd BAB,
$$
from which the result follows (with $c=(n\al)\inv\be)$.\qed

If, in the context of this lemma, we set $B$ equal to $A\snv$, then
$\nab$ is the Bose-Mesner algebra of an association scheme.  So the
assumption $A\in\nab$ implies that the diagonal of $A$ is constant and
the column sums of $B$ are equal.  Therefore, \ref{nab-jpr} extends
Proposition~9 of Jaeger, Matsumoto and Nomura \cite{MR99f:05125}.
(This proposition asserts that $cA$ is a spin model for some nonzero
$c$ if and only if $A\in\na$.)

\section{gauges}
Gauge Equivalence

We show here that each member of an invertible Jones pair almost
determines the other.

If $D$ is an invertible diagonal matrix, we call $D\inv JD$ a {\sl
dual permutation matrix}.  We note that if $A\snv\circ C=D\inv JD$,
then
$$
C =A\circ(D\inv JD) =D\inv AD.
$$
Thus $A\snv\circ C$ is a dual permutation matrix if and only if $A$
and $C$ are diagonally equivalent.  The Schur inverse of a dual
permutation matrix is a dual permutation matrix.  (The concept of dual
permutation matrix comes from Jaeger and Nomura \cite{MR1723186}.)

\lemma{dlperm}
If $A$, $C$ and $M$ are Schur invertible matrices and $\xa\del M =\del
MX_C$, then $C\snv\circ A$ is a dual permutation matrix.  If $B$, $C$ and
$M$ are invertible matrices and $\del BX_M=X_M\del C$, then $CB\inv$ is
a permutation matrix.

\proof
Suppose $\xa\del M =\del MX_C$.  Applying each side of this equality to
$\eij$ yields
$$
M_{i,j}\,Ae_i =Ce_i\circ Me_j
$$
and, as $C\snv$ exists, we get
$$
Me_j =M_{i,j} (C\snv\circ A)e_i.\eqdef{mmca}
$$
Therefore, each column of $(C\snv\circ A)$ is a multiple of $Me_1$,
and so $C\snv\circ A =D_1JD_2$, for some invertible diagonal matrices
$D_1$ and $D_2$.  Next, Equation~\eqref{mmca} implies that
$$
M_{i,j}=e_i^TMe_j=M_{i,j}\,e_i^T(C\snv\circ A)e_i
$$
and, as $M_{i,j}\ne0$ this implies that $e_i^T(C\snv\circ A)e_i=1$.
Hence $D_1D_2=I$ and $C\snv\circ A$ is a dual permutation matrix.

Now, suppose $\del BX_M=X_M\del C$.  Then
we get
$$
Me_i\circ Be_j =C_{i,j}\,Me_i
$$
and so, if we multiply both sides of this by $(B\inv)_{j,k}$ and sum
over $j$, we get
$$
M_{k,i}e_k=Me_i\circ e_k =\left(\sum_k C_{i,j}(B\inv)_{j,k}\right)Me_i
	=(CB\inv)_{i,k}\,Me_i.\eqdef{mbcm}
$$
This implies that each column of $M$ is a multiple of some vector
$e_r$.  Since $M$ is invertible, no column is zero, and so we also see
that, for value of $i$ there is at most one index $k$ such that
$(CB\inv)_{i,k}\ne0$.  If $M_{k,i}\ne0$, then \eqref{mbcm} gives
$$
M_{k,i} =(CB\inv)_{i,k}\,e_k^TMe_i =(CB\inv)_{i,k}\,M_{k,i}
$$
and thus $(CB\inv)_{i,k}=1$.  Since $CB\inv$ is invertible, we conclude
it is a permutation matrix.\qed

\corollary{}
Let $(A,B)$ be an invertible one-sided Jones pair.  If $(C,B)$ is also
an invertible one-sided Jones pair there is an invertible diagonal
matrix $D$ such that $C=D\inv AD$.

\proof
By \ref{a3b3}, if $(A,B)$ and $(C,B)$ are invertible one-sided Jones
pairs, then
$$
\xdx AAB =\dxd AB{B^T},\qquad \xdx CCB =\dxd CB{B^T}
$$
and consequently
$$
\xa\del{A}\del{C\snv} X_{C\inv}=\del{A}\del{C\snv}.
$$
Therefore,
$$
\xa\del{C\snv\circ A}=\del{C\snv\circ A}X_C
$$
and, by \ref{dlperm}, this yields that $C\snv\circ A$ is a dual
permutation matrix.  Hence $A$ and $C$ are diagonally equivalent.\qed

\corollary{}
Let $(A,B)$ be an invertible one-sided Jones pair.  If $(A,C)$ is also
an invertible one-sided Jones pair there is a permutation matrix $P$
such that $C=BP$.

\proof
Now we have 
$$
\xdx AAB =\dxd AB{B^T},\qquad \xdx AAC =\dxd AC{C^T},
$$
and consequently
$$
X_{C\inv}X_B=\del{C\snt} X_{C\inv}X_B\del{B^T}.
$$
Hence
$$
\del{C^T}X_{C\inv B}=X_{C\inv B}\del{B^T}
$$
and \ref{dlperm} implies that $B^TC^{-T}$ is a permutation
matrix.\qed

These results have interesting consequences.  We saw that $(A^T,B)$ is
a one-sided Jones pair if $(A,B)$ is, so we may deduce that if $(A,B)$
is an invertible one-sided Jones pair, then there is a diagonal matrix
$C$ such that $C\inv AC= A^T$.  Further, if we are working over $\cx$
and $C\inv AC=A^T$ then there is a diagonal matrix $C_1$ such that
$C_1^2=C$; then we have
$$
C_1\inv AC_1=C_1A^TC_1\inv=(C_1\inv AC_1)^T
$$
and so $C_1\inv AC_1$ is symmetric.

If $(A,B)$ is an invertible Jones pair, then so is $(A,B^T)$, whence it
follows that $B^T=BP$, for some permutation matrix $P$.  Since
$$
B=(B^T)^T =(BP)^T =P^TB^T =P^TBP
$$
we see that $P$ must commute with $B$.  If $P$ has odd order, then
there is an integer $r$ such that $P^{2r}=P$.  Suppose $Q=P^r$.  Then
$Q$ commutes with $B$ and $Q^TB^T=BQ$.  Therefore, $BQ$ is symmetric.

(The facts that if $(A,B)$ is an invertible Jones pair, then
$A^T=C\inv AC$ for some diagonal matrix $C$ and $B^T=BP$ for some
permutation matrix $P$ are due to Jaeger \cite{MR1771614}, but our
proofs are new, and simpler.)

\section{4wt}
Four-Weight Spin Models

A {\sl four-weight spin model}\/ is a 5-tuple $(W_1,W_2,W_3,W_4;D)$
where $W_1$, $W_2$, $W_3$, $W_4$ are $n\times n$ complex matrices and
$D$ is a square root of $n$ such that:
$$\eqalignno{
	W_1\circ W_3^T=J,&\quad W_2\circ W_4^T=J,&\eqdef{4schur}\cr
	W_1W_3=nI,&\quad W_2W_4=nI,&\eqdef{4inv}\cr
	\sum_{h=1}^n(W_1)_{k,h}(W_1)_{h,i}(W_4)_{h,j} 
		&=D(W_4)_{i,j}(W_1)_{k,i}(W_4)_{k,j},&\eqdef{type3a}\cr
	\sum_{h=1}^n(W_1)_{h,k}(W_1)_{i,h}(W_4)_{j,h} 
		&=D(W_4)_{j,i}(W_1)_{i,k}(W_4)_{j,k}.&\eqdef{type3b}\cr}
$$
Note that \eqref{type3a} and \eqref{type3b} are equivalent to (3a) and
(3b) in \cite[p.~1]{MR97f:57005} respectively.  The original spin
models due to Jones \cite{MR89m:57005} are referred to as {\sl
two-weight}\/ spin models.  

\theorem{}
Suppose that $A, B\in M_n(\cx)$ and $D^2=n$.  Then the following are
equivalent. 
\item{(i)} $(A,B)$ is an invertible Jones pair.
\item{(ii)} $(DA,nB\inv, DA\inv, B;D)$ is a four-weight spin model.

\proof
Write $A=D\inv W_1$ and $B=W_4$.  Then \eqref{type3a} is equivalent to
\ref{evecs}, and holds if and only if $(A,B)$ is a one-sided Jones
pair.  Observe next that \eqref{type3b} is obtain from \eqref{type3a}
by replacing $W_1$ and $W_4$ by their transposes.  This implies that
\eqref{type3b} holds if and only if $(A^T,B^T)$ is a one-sided Jones 
pair which, by \ref{jp-equivs}, is equivalent to $(A,B^T)$ being a
one-sided Jones pair.  This shows that (ii) implies (i).  The converse
is easy, given \ref{jp-type2}.\qed

Thus invertible Jones pairs are equivalent to four-weight spin models.
Our treatment shows that most of the theory developed in
\cite{bet-bm4wt,MR97f:57005} holds under the weaker
assumption that $(A,B)$ is a one-sided Jones pair.

Four-weight spin models were introduced as a generalization of the
generalized spin models of Kawagoe, Munemasa and Watatani
\cite{MR96g:57009}.  In our terms, Kawagoe et al defined a generalized
spin model to be an $n\times n$ type II matrix $A$ such that
$(A^T,\sqrt{n}A\snv)$ is a one-sided Jones pair.  We use the exchange
lemma to show that such a pair must be two-sided.  (This shows that
four-weight spin models are indeed a generalization of generalized
spin models, as noted in Bannai and Bannai \cite{MR97f:57005}.)

\lemma{}
If $A$ is type II and $(A^T,\sqrt{n}A\snv)$ is a one-sided Jones pair,
then it is a Jones pair.

\proof
By \ref{a3b3},
$$
\xdx {A^T}{A^T}{A\snv} =\dxd{A^T}{A\snv}{\sqrt{n}A\snt}.
$$
Inverting each side of this, we find that
$$
\xdx{A^{(-)-1}}{A\snt}{A^{-T}}
	={1\over\sqrt n}\dxd{A^T}{A^{(-)-1}}{A\snt}
$$
Since $A^{-T}={1\over n}A\snv$, we get
$$
\xdx{A^T}{A\snt}{A^{-T}} ={1\over\sqrt n}\dxd{A^T}{A^T}{A\snt}
$$
which gives
$$
\dxd{\sqrt{n}A\snt}{A^T}{\sqrt{n}A\snt} 
	=\xdx{A^T}{\sqrt{n}A\snt}{A^T}.
$$
This implies that $(A^T,\sqrt{n}A\snt)$ is a one-sided Jones pair.\qed

\section{evecs2}
Algebras and Bijections

Given a one-sided Jones pair we have a number of algebras including
$\nab$, $\nabp$ and the Bose-Mesner algebra $\na$.  In this section we
study some of the relations between these.

\theorem{ABC}
Let $A$, $B$ and $C$ be type II matrices with the same order.
If $F\in\nab$ and $G\in\nom{B\snv,C}$, then $F\circ G\in\nom{A,C}$ and
$$
\Th_{A,C}(F\circ G)=n\inv\thab(F)\,\Th_{B\snv,C}(G).
$$

\proof
If $\xdx FBA=\dxd BA{F'}$, then, by \ref{rbabas}(b),
$$
\dxd F{B\snv}{B^T}=\xdx A{A\inv}{F'}.\eqdef{fbb}
$$
Similarly, applying \ref{rbabas}(b) to $\xdx GC{B\snv}
=\dxd{C}{B\snv}{G'}$,  we get
$$
\dxd G{C\snv}{C^T} =\xdx{B\snv}{B\snv{}\inv}{G'} 
	=n\inv\xdx{B\snv}{B^T}{G'}.\eqdef{gcc}
$$
>From \eqref{fbb}, we have
$$
\del{F}=\xdx A{A\inv}{F'}(X_{B\snv}\del{B^T})\inv,
$$
which, combined with \eqref{gcc}, yields
$$
\del{F}\dxd G{C\snv}{C^T} =n\inv\xdx A{A\inv}{F'}X_{G'},
$$
and therefore
$$
\dxd{F\circ G}{C\snv}{C^T} =n\inv\xdx A{A\inv}{F'G'}.
$$
By applying the exchange lemma to this,
we find that
$$
\dxd{C\snv}{F\circ G}C=n\inv\xdx{A}{F'G'}{A\inv},
$$
whence
$$
\xdx{F\circ G}CA=n\inv\dxd CA{F'G'}.
$$
Therefore, $F\circ G\in\nom{A,C}$ and 
$$
\Th_{A,C}(F\circ G)=n\inv\thab(F)\,\Th_{B\snv,C}(G).\qed
$$

\lemma{nabt}
Suppose $A$ and $B$ are type II matrices of the same order.  If $G\in\nab$,
then $G^T\in\nom{A\snv,B\snv}$ and
$\Th_{A\snv,B\snv}(G^T)=\thab(G)$.

\proof
Suppose $R\in\nab$ and $\xdx RBA=\dxd BAS$.  Then
$$
\xdx{A\inv}{B\snv}R =\dxd S{A\inv}{B\snv}
$$
and, taking the  transpose of this, we find that
$$
\xdx{R^T}{B\snv}{A^{-T}}=\dxd{B\snv}{A^{-T}}S.
$$
As $A$ is type II we have $nA^{-T}=A\snv$, and we conclude that
$R^T\in\nom{A\snv,B\snv}$ and $\Th_{A\snv,B\snv}(R^T)=S$.\qed

As an example, suppose $A$ is type II and $G\in\na$.  Then by the
lemma, 
$$
G^T \in\nom{A\snv,A} =\na
$$ 
and
$$
\Th_{A,A\snv}(G) =\Th_{A\snv,A}(G^T) =\Th_{A,A\snv}(G^T)^T,
$$
which implies that $\Th_A(G^T)=\Th_A(G)^T$.

\theorem{nanab}
Suppose $A$ and $B$ are type II matrices of the same order.  If $F\in\na$,
$G\in\nab$ and $H\in\nom{B}$, then $F\circ G$ and $G\circ H$ lie in
$\nab$ and
$$\eqalign{
\thab(F\circ G)&=n\inv\Th_A(F)\,\thab(G),\cr
\thab(G\circ H)&=n\inv\thab(G)\,\Th_{B}(H)^T,\cr}
$$

\proof
If we apply \ref{ABC} to the triple $(A,A\snv,B)$, we get the first
claim. 

For the second claim, note that if $H\in\nom{B}$, then
$H\in\nom{B\snv}$ and $\Th_{B\snv}(H)=\Th_B(H)^T$.  Now, if we apply
\ref{ABC} to the triple $(A,B,B)$, then we have
$$
\thab(G\circ H) =n\inv\thab(G)\,\Th_{B\snv}(H),
$$
from which the second assertion follows.\qed

This corollary shows that $\na\circ\nab\sbs\nab$ and
$\nom{A^T}\nabp\sbs\nabp$.  We will see below that equality holds.

\theorem{ght}
Suppose $A$ and $B$ are type II matrices of the same order.  If $F$
and $G$ lie in $\nab$, then $F\circ G^T\in\na\cap\nom{B}$ and
$$\eqalign{
\Th_A(F\circ G^T) &=n\inv\thab(F)\,\thab(G)^T,\cr
\Th_B(F\circ G^T) &=n\inv \thab(F)^T\,\thab(G).\cr}
$$
 
\proof
If $G\in\nab$, then, by \ref{nabt}, $G^T\in\nom{A\snv,B\snv}$ and
$$
\Th_{A\snv,B\snv}(G^T)=\thab(G).
$$
If we now apply \ref{ABC} to the triple $(A,B,A\snv)$, we find that
$F\circ G^T\in\na$ and that $\Th_A(F\circ G^T)$ is as stated.
For the remaining claims, apply \ref{ABC} with the
triple $(B,A,B\snv)$.\qed

Our next result is an easy consequence of \ref{ght}.  It implies
that if $(A,B)$ is an invertible one-sided Jones pair and $A$ is
symmetric, then either $A\circ A$ is a linear combination of $I$ and
$J$, or $\na$ is the Bose-Mesner algebra of an association scheme with
at least two classes.  (We will discuss the first case again in the
final section.)

\corollary{acat-bbt}
If $(A,B)$ is an invertible one-sided Jones pair, then $A\circ
A^T\in\na$ and $\Th_A(A\circ A^T) =n\inv BB^T$.\qed

Among other important consequences, our next result implies that if
$(A,B)$ is an invertible Jones pair, then $\na$ and $\nab$ have the
same dimension.

\theorem{aona}
Suppose $A$ and $B$ are type II matrices of the same order.
If $\nab$ contains a Schur invertible matrix $G$
then $\na=\nom{B}$,
\smallbreak
\item{(a)} $G\circ\na=\nab$ and $ G^T\circ\nab=\na$.
\item{(b)} If $H=\thab(G)$, then $\nom{A^T}H=\nabp$ and 
$\nabp H^T=\nom{A^T}$.

\proof
By \ref{nanab} we see that $\na\circ G\sbs\nab$, while \ref{ght}
implies that $\nab\circ G^T\sbs\na$.  Since $G\snv$ exists, Schur
multiplication by $G$ is injective, and so (a) follows.

\ref{nanab} also implies that $G\circ\nom{B}\sbs\nab$, while
\ref{ght} implies that $G^T\circ\nab\sbs\nom{B}$.  Hence 
$G\circ\nom{B}=\nab$ and $G^T\circ\nab=\nom{B}$, and therefore
$\na=\nom{B}$. 

By \ref{nanab}
$$
\thab(G\circ\na) =n\inv \Th_A(\na)H.
$$
By \ref{nancor}, $\Th_A(\na)=\nom{A^T}$ and so the first part of the
second claim follows.  From \ref{ght} we see that
$$
\Th_A(\nab\circ G^T) =n\inv\nabp H^T,
$$
whence the second part of the second claim follows.\qed

If $(A,B)$ is an invertible one-sided Jones pair, then $A$ is a Schur
invertible element of $\nab$.  Thus we have immediately:

\corollary{}
If $(A,B)$ is an invertible one-sided Jones pair, then
$A^T\circ\nab=\na$ and $\nab=A\circ\na$.\qed

Since $\na$ is a Bose-Mesner algebra, this result implies that all
matrices in $\nab$ have constant diagonal.  Since $\na$ is closed
under transposes. we also see that if $A$ is symmetric, then $\nab$ is
closed under transposes and (by \ref{nabt}) that $\nab=\nom{A\snv,B\snv}$.

\corollary{}
If $(A,B)$ is an invertible one-sided Jones pair, then
$\dim\na=\dim\nab$ and $\nom{A^T}=\na=\nom{B}$.

\proof
Given our hypothesis, $A$ and $A^T$ are Schur invertible, and it
follows from the first part of \ref{aona} that $\dim\na=\dim\nab$ and
that $\na=\nom{B}$.  Since $(A^T,B)$ is an invertible one-sided Jones
pair if $(A,B)$ is, we also have $\nom{A^T}=\nom{B}$.\qed

Etsuko Bannai\cite{bet-bm4wt} proved that
$\na=\nom{B}=\nom{A^T}=\nom{B^T}$ for four-weight spin models; this
extended unpublished work by H. Guo and T. Huang, who had shown that
$\na=\nom{A^T}$.  Note that if $(A,B)$ is an invertible Jones pair,
then $(A,B^T)$ is an invertible Jones pair, and so our previous result
also implies that $\na=\nom{B^T}$ in this case.

\corollary{thathb}
If $(A,B)$ is an invertible Jones pair and $F\in\na=\nom{B}$,
then $\Th_B(F)^T =B\inv\Th_A(F)B$.

\proof
We have seen that $\na=\nom{B}$ and therefore $A\circ\nom{B}=\nab$.
Suppose $F\in\nom{B}$.  Then by \ref{nanab},
$$
\thab(A\circ F)=n\inv B\Th_B(F)^T
$$
and
$$
\thab(A\circ F)=n\inv\Th_A(F)B.
$$
Hence the result follows.\qed

This result implies that $B\inv\na B=\na$.  Similarly, if $(A,B)$ is
an invertible Jones pair, then we also find that $B^{-T}\na B^T=\na$.
In this case $B^T=BP$ for some permutation matrix $P$, and $P^T\na
P=\na$.

\section{2ass}
A Dual Pair of Schemes

Let $A$ and $B$ be type II matrices of the same order, and let $W$ be
the matrix defined by
$$
W:=\pmatrix{A&B\snv\cr -A&B\snv\cr}.
$$
Then it is easy to verify that $W$ is a type II matrix; we are going
to describe its Nomura algebra.

\theorem{}
If $A$ and $B$ are type II matrices of the same order and
$$
W=\pmatrix{A&B\snv\cr -A&B\snv\cr},
$$
then $\nom{W}$ consists of the matrices
$$
\pmatrix{F+R&F-R\cr F-R&F+R},
$$
where $F\in\na\cap\nom{B}$ and $R\in\nom{A,B}\cap\nom{A\snv,B\snv}$.

\proof
Suppose $M$, $N$, $P$ and $Q$ are $n\times n$ matrices.  Then
the matrix
$$
Z:= \pmatrix{M&N\cr P&Q\cr}
$$
lies in $\nom{W}$ if and only if all of the following vectors are
eigenvectors for $Z$:
$$\displaylines{
\pmatrix{Ae_i\circ A\snv e_j\cr Ae_i\circ A\snv e_j\cr},\quad
\pmatrix{B\snv e_i\circ Be_j\cr B\snv e_i\circ Be_j\cr},\cr
\pmatrix{Ae_i\circ Be_j\cr -Ae_i\circ Be_j\cr},\quad
\pmatrix{A\snv e_i\circ B\snv e_j\cr -A\snv e_i\circ B\snv e_j\cr}}
$$
The first two of these four sets of vectors are eigenvectors for $Z$ if
and only if both $M+N$ and $P+Q$ lie in $\na$ and 
$$
\Th_A(M+N)=\Th_A(P+Q).
$$
This last condition implies that $M+N=P+Q$; we may assume both sums
equal $F$, where $F\in\na$.  Similarly the second set of vectors
consists of eigenvectors for $Z$ if and only if $M+N$ and $P+Q$ lie in
$\nom{B}$, and therefore $F\in\na\cap\nom{B}$.

The third and fourth sets of vectors are eigenvectors for $Z$ if and
only if $M-N$ and $P-Q$ lie in $\nom{A,B}\cap\nom{A\snv,B\snv}$ and 
$$
\thab(M-N)=\thab(P-Q),
$$
whence $M-N=P-Q=R$.\qed

This result shows that $\nom{W}$ is the direct sum of two subspaces.
The first consists of the matrices of the form
$$
\pmatrix{F&F\cr F&F\cr}\eqdef{ffff}
$$
where $F\in\na\cap\nom{B}$.  This set of matrices is closed under
multiplication and Schur multiplication (but does not contain $I$).
The second subspace consists of the matrices
$$
\pmatrix{R&-R\cr -R&R\cr}\eqdef{rrrr}
$$
where $R\in\nom{A,B}\cap\nom{A\snv,B\snv}$.

\corollary{}
Suppose $F\in\na\cap\nom{B}$ and $R\in\nom{A,B}\cap\nom{A\snv,B\snv}$.
If 
$$
Z :={1\over2}\pmatrix{F+R&F-R\cr F-R&F+R\cr}
$$
then
$$
\Th_W(Z)=\pmatrix{
	\Th_A(F)&\thab(R)\cr
	\Th_{B\snv,A\snv}(R)&\Th_{B\snv}(F)\cr}.
$$

\proof
Assume $A$ and $B$ are $n\times n$ matrices and $1\le i,j\le n$.
Then, for example
$$
We_{i+n}\circ W\snv e_j =\pmatrix{
	B\snv e_i\circ A\snv e_j\cr
	-B\snv e_i\circ A\snv e_j\cr}.
$$
This determines the $(2,1)$-block of $\Th_W(Z)$, and the other blocks
can be found in a similar way.\qed

Observe that $\Th_{B\snv}(F)=(\Th_B(F))^T$ and so, by
\ref{thathb},
$$
\Th_{B\snv}(F)=B\inv\Th_A(F)B.
$$
By \ref{nabt} we have that $\nom{A\snv,B\snv}=\nab^T$
and
$$
\Th_{A\snv,B\snv}(R)=\thab(R^T).
$$

We can now state one of the main conclusions of our paper.  

\corollary{}
Let $(A,B)$ be an invertible Jones pair of $n\times n$ matrices.
Assume $A$ is symmetric and $\dim\na=m$.  Let $W$ be the $2n\times2n$
type II matrix defined above.  Then $\nom{W}$ is the Bose-Mesner
algebra of an imprimitive association scheme with $2m-1$ classes
which contains the matrix
$$
\pmatrix{A&-A\cr -A&A\cr}.
$$
The image of this under $\Th_W$ is the following matrix in the dual
scheme (with Bose-Mesner algebra $\nom{W^T}$):
$$
\pmatrix{0&B\cr B^T&0}.
$$

\proof
Since $(A,B)$ is an invertible Jones pair, $\na=\nom{B}$.  Hence
the space of matrices in \eqref{ffff} has dimension $m$.  Since $A$ is
symmetric, $\nab$ is closed under transposes and therefore equals
$\nom{A\snv,B\snv}$.  Hence the space of matrices in \eqref{rrrr} also
has dimension $m$, and consequently $\dim\nom{W}=2m$.

Because $A\in\nab$, it follows that
$$
\pmatrix{A&-A\cr -A&A\cr}\in\nom{W}
$$
and then the previous corollary implies that
$$
\pmatrix{0&B\cr B^T&0}\in\nom{W^T}.\qed
$$

\section{dim2}
Dimension Two

We define the {\sl dimension}\/ of a Jones pair $(A,B)$ to be the
dimension of $\nab$ and we define the {\sl degree}\/ to be the number
of distinct entries in $B$ or, equivalently, the number of distinct
eigenvalues of $A$.  Since $A\in\nab$, the degree is bounded above by
the dimension.  Since $\nab$ contains $I$ and $A$ the dimension of a
pair is at least two, unless $A=I$ (and then $B=J$).

Suppose $(A,B)$ is an invertible Jones pair of $n\times n$ matrices
with dimension two, and that $A$ is symmetric.  Since
$\dim\na=\dim\nab$, we see that $\na$ is the span of $I$ and $J$, so
$\na$ is the Bose-Mesner algebra of an association scheme with one
class It follows that there are complex numbers $a$ and $b$ such that
$$
A\circ A =aI+bJ.
$$
Therefore there is symmetric matrix $C$ such that $C\circ I=0$ and
$C_{i,j}=\pm1$ if $i\ne j$ and $A=\la I+\ga C$.  Further, since $\nab$
has dimension two, the minimal polynomial of $A$ is quadratic.  Hence
the minimal polynomial of $C$ is quadratic and, since $(C^2)\circ
I=(n-1)I$, there is an integer $\de$ such that
$$
C^2-\de C-(n-1)I=0.
$$
This implies that $C$ is the matrix of a regular two-graph.  (For more
information on regular two-graphs, see Seidel's two surveys in
\cite{MR92m:01098}, and for more on the connection with type II
matrices, see \cite{cg-sd2g}.)

If $A$ has quadratic minimal polynomial, then $B$ has exactly two
distinct entries and so is a linear combination of $J$ and a
$01$-matrix $N$.  In \cite{bansaw-4wtdes}, Bannai and Sawano show that
$N$ must be the incidence matrix of a symmetric design, and
characterize the designs that can arise in this way.  It is well known
that symmetric designs on $n$ points correspond to bipartite distance
regular graphs on $2n$ vertices with diameter three, and less well
known that a formal dual of such a scheme is the association scheme
associated to a regular two-graph.

Finally, if $(A,B)$ is a one-sided Jones pair and $A$ has quadratic
minimal polynomial then the algebra generated by $\xa$ and $\deb$ is a
quotient of the Hecke algebra.  It follows that the link invariant we
obtain is a specialization of the homfly polynomial.  (For this see
Jones \cite[Section~4]{MR89c:46092}.)

\nonumsection%
References

\bibliography{joprs4}
\bibliographystyle{acm}

\bye